\documentclass[11pt]{amsart}

\usepackage{amsfonts,latexsym,graphicx,fullpage}

\newcommand\R{\mathbb{R}}

\newcommand\G{S}
\def\d{\partial}

\newtheorem{proposition}{Proposition}

\newtheorem{lemma}[proposition]{Lemma}
\newtheorem{corollary}[proposition]{Corollary}

\theoremstyle{definition}
\newtheorem{definition}[proposition]{Definition}

\theoremstyle{remark}
\newtheorem{remark}[proposition]{Remark}

\numberwithin{proposition}{section}

\theoremstyle{plain}

\begin{document}

\title{The Legendrian Satellite Construction}
\author{Lenhard L.\ Ng}
\address{American Institute of Mathematics, 360 Portage Avenue,
Palo Alto, CA 94306}
\curraddr{School of Mathematics, Institute for Advanced Study,
Einstein Drive, Princeton, NJ 08540}
\email{ng@ias.edu}
\urladdr{http://www.math.ias.edu/\~{}ng/}
\date{December 2001}
\thanks{The author is supported by an American
Institute of Mathematics Five-Year Fellowship}


\begin{abstract}
We examine the Legendrian analogue of the topological satellite
construction for knots, and deduce some results for specific Legendrian
knots and links in standard contact three-space and the solid torus.  
In particular, we show that the
Chekanov-Eliashberg contact homology invariants of Legendrian Whitehead
doubles of stabilized knots contain no nonclassical information.
\end{abstract}

\maketitle

\vspace{24pt}

\section{Introduction}

The spaces $\R^3$ and $S^1\times\R^2$ both have a standard contact
structure given by the kernel of the $1$-form $dz - y\,dx$, where
we view the solid torus $S^1\times\R^2$ as $\R^3$ modulo the relation
$(x,y,z) \sim (x+1,y,z)$.  We will assume that the reader is familiar
with some basic concepts in Legendrian knot theory, such as front projections,
the Thurston-Bennequin number, and the rotation number;
see, e.g., \cite{bib:Ng}, which we will use extensively.

The problem of classifying Legendrian knots in standard contact $\R^3$
up to Legendrian isotopy has attracted much recent attention.
In this note, we study one particular construction on Legendrian knots,
the Legendrian satellite construction, which relates knots in
$\R^3$ and in $S^1\times\R^2$.  This is the Legendrian analogue of
the satellite construction in the smooth category, which
glues a link in the solid torus $S^1\times\R^2$ into a tubular neighborhood
of a knot in $\R^3$ to produce a link in $\R^3$.
We examine some consequences of Legendrian satellites for
the Legendrian knot classification problem in both $\R^3$ and 
$S^1\times\R^2$; in particular, we recover previously-known
results for knots in both spaces, and prove a new result
in $S^1\times\R^2$ (Proposition~\ref{prop:whiteheadknot}).

\begin{figure}[t]
\centering{
\includegraphics[width=0.7in,angle=270]{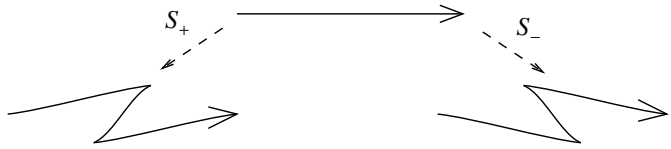}
}
\caption{Stabilization of a Legendrian link, in the front projection.}
\label{fig:stabilization}
\end{figure}

The motivation for this work is that Legendrian satellites may
provide nontrivial, nonclassical invariants of stabilized
Legendrian knots in $\R^3$.  Here we recall that there are two stabilization
operators $S_\pm$ on Legendrian knots, decreasing $tb$ by $1$ and
changing $r$ by $\pm 1$, which replace a segment of the knot's front
projection by a zigzag, as shown in Figure~\ref{fig:stabilization}.
Understanding stabilized knots is an important open problem in
Legendrian knot theory; it also has repercussions for the classification
of transverse knots.

It seems possible that satellites of stabilized knots may contain
interesting information through the Chekanov-Eliashberg differential
graded algebra invariant \cite{bib:Che,bib:Ng}, which is derived
from contact homology \cite{bib:Eli}.  We will show that,
unfortunately, the DGAs of the simplest Legendrian satellites of
stabilized knots do not encode any useful information.  The computation
used in the proof may be of interest as the first involved computation
manipulating the DGA invariant directly, rather than using easier invariants
such as Poincar\'e polynomials \cite{bib:Che}
or the characteristic algebra \cite{bib:Ng}.  In any case, more complicated
satellites may well give nonclassical invariants of stabilized knots,
as has been suggested by Michatchev \cite{bib:Mi}.

We define the construction in Section~\ref{sec:gluing}, and show
how it immediately implies facts about solid-torus links, including
some that could not be shown using any previously known techniques.
In Section~\ref{sec:doubles}, we show the result mentioned above
about DGAs of satellites of stabilized knots;
the key step is Lemma~\ref{lem:doubleproof}, which is proven in
Section~\ref{sec:prooflemdouble}.

\vspace{11pt}

\noindent
{\it Acknowledgments.}
I would like to thank Yasha Eliashberg, John Etnyre, 
Kirill Michatchev,
Tom Mrowka, Josh Sabloff, and Lisa Traynor for helpful discussions, and
the American Institute of Mathematics
for sponsoring the fall 2000 Low-Dimensional Contact Geometry program, 
during which most of these discussions took place.
This note was originally part of my MIT Ph.D.\ thesis, which
was written with the support of research assistantships from grants
from the National Science Foundation and Department of Defense.

\section{Construction}
\label{sec:gluing}

\begin{figure}
\centering{\includegraphics[width=5in,angle=270]{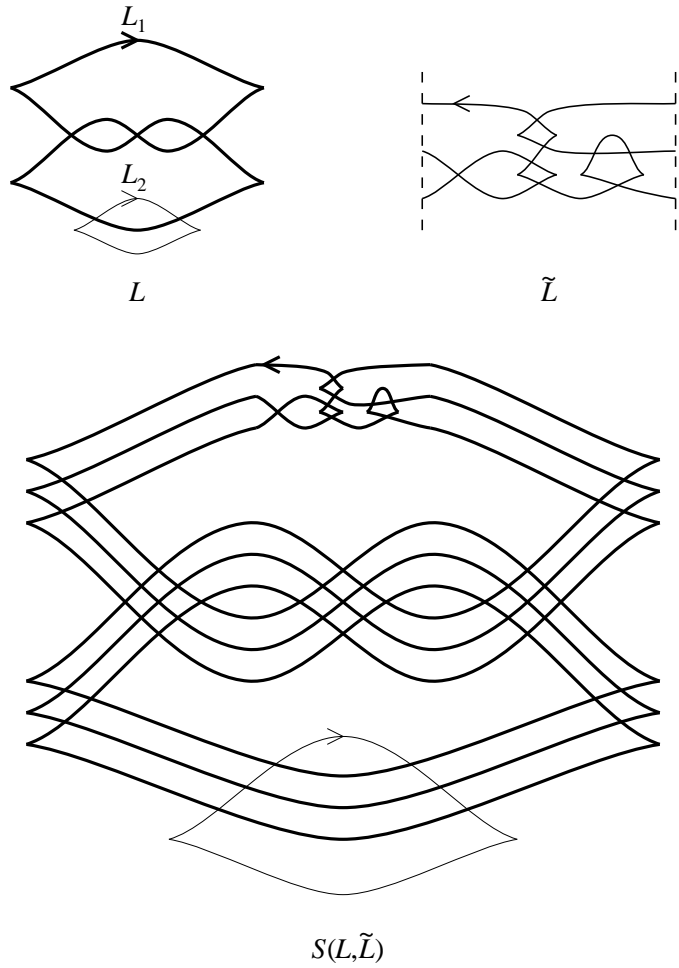}}
\caption{
Gluing a solid-torus link $\tilde{L}$ into an $\R^3$ link $L$,
to form the satellite link $\G(L,\tilde{L})$.
}
\label{fig:torusgluing}
\end{figure}

As in $\R^3$, Legendrian links
in the solid torus $S^1\times\R^2$ may be represented by their front 
projections to the
$xz$ plane, with the understanding that the $x$ direction is now periodic.
If we view $S^1 \times \R^2$ as $[0,1] \times \R^2$ with $\{0\} \times \R^2$
identified with $\{1\} \times \R^2$, then we can draw the front projection
of a solid-torus Legendrian link as a front in $[0,1] \times \R$
with the two boundary components identified.  We depict the boundary
components by dashed lines;
see Figure~\ref{fig:torusgluing} for an illustration.
For a Legendrian link $\tilde{L}$ in the solid torus, let the
{\it endpoints} of $\tilde{L}$ be $\tilde{L} \cap (\{0\} \times \R^2)$,
that is, the points where the front for $\tilde{L}$ intersects the dashed
lines.

\begin{remark} {\it Invariants of solid-torus links.}
There are three classical invariants of links on the solid torus:
the Thurston-Bennequin number $tb$ and rotation number $r$ can be calculated
from the front of a solid-torus link exactly as in $\R^3$; and the winding
number $w$ is the number of times the link winds around the $S^1$
direction of $S^1 \times \R^2$.  Clearly the $tb$, $r$, and $w$ associated
to any subset of the components of a solid-torus link also 
give invariants of the link.

In \cite{bib:NT}, L.\ Traynor and the author show that the 
Chekanov-Eliashberg DGA
can be defined for links on the solid torus, thus yielding a nonclassical
invariant.  For certain links with two components, \cite{bib:Tra}
defines another nonclassical invariant based on generating functions.
We will give examples in this section of solid-torus knots which
are not Legendrian isotopic, but which cannot be distinguished
using any of these invariants.
\end{remark}

We now introduce the Legendrian satellite construction.
Let $L$ be an oriented Legendrian link in $\R^3$ with one distinguished
component $L_1$, and let $\tilde{L}$ be an oriented Legendrian link
in $S^1 \times \R^2$.  We give two definitions of the Legendrian 
satellite $\G(L,\tilde{L}) \subset \R^3$, one abstract, one concrete.

A tubular neighborhood of $L_1$ is a solid torus;
the characteristic foliation on the boundary of this torus wraps
around the torus $tb(L_1)$ times.  By cutting the tubular neighborhood
at a cross-sectional disk, untwisting it $tb(L_1)$ times, and 
regluing, we obtain a solid torus contactomorphic to
$S^1 \times \R^2$ with the standard contact structure.
Thus we can embed $\tilde{L} \subset S^1 \times \R^2$ as a Legendrian 
link in a tubular neighborhood of $L_1$.  Replacing the component
$L_1$ in $L$ by this new link gives $\G(L,\tilde{L})$.

We can redefine $\G(L,\tilde{L})$ in terms of the fronts for $L$ and
$\tilde{L}$.  First, we recall the definition of an $n$-copy from
\cite{bib:Mi}.

\begin{definition}
Given a Legendrian knot $K$, its {\it $n$-copy} is the link consisting
of $n$ copies of $K$ which differ from each other through small perturbations
in the transversal direction.  In the front projection, the $n$-copy consists
of $n$ copies of $K$, differing from each other by small shifts in the $z$
direction.  The $2$-copy is also known as the {\it double}.
\end{definition}
 
Now suppose that $\tilde{L}$ has $n$ endpoints.  By cutting along the dotted 
lines (i.e., the endpoints of $\tilde{L}$), we can
embed $\tilde{L}$ as a Legendrian tangle in $\R^3$ with $2n$ ends.
Replace the front of the first component $L_1$ of $L$ by the $n$-copy of 
$L_1$.  Then
choose a small segment of $L_1$ which is oriented from left to right;
excise the corresponding $n$ pieces of the $n$-copy of $L_1$, and
replace them by the front for $\tilde{L}$, cut along its endpoints.
See Figure~\ref{fig:torusgluing} for an illustration.

\begin{definition}
The resulting link $\G(L,\tilde{L}) \subset \R^3$ is the 
{\it Legendrian satellite}
of $L \subset \R^3$ and $\tilde{L} \subset S^1\times\R^2$.  
We give $\G(L,\tilde{L})$
the orientation derived from the orientations on $\tilde{L}$ (for
the glued $n$-copy of $L_1$) and on $L$ (for the components of
$L$ besides $L_1$).
\end{definition}

The Legendrian satellite construction is motivated by the special 
case of Whitehead doubles (see Section~\ref{sec:doubles}),
which were introduced by Eliashberg and subsequently used by
Fuchs \cite{bib:Fuc}.

\begin{remark} {\it Classical invariants of Legendrian satellites.}
Before we show that $\G(L,\tilde{L})$ is well-defined up to Legendrian 
isotopy, we note that the classical invariants of
$\G(L,\tilde{L})$ are easily computable from those of $L$ and $\tilde{L}$.
Indeed, a straightforward computation with front diagrams yields
\begin{gather*}
tb(\G(L,\tilde{L})) = (w(\tilde{L}))^2 \,tb(L) + tb(\tilde{L}) \\
r(\G(L,\tilde{L})) = w(\tilde{L}) r(L) + r(\tilde{L})
\end{gather*}
when $L$ is a knot, with a similar but slightly more complicated
formula when $L$ is a multi-component link.
\label{rmk:satellitetbr}
\end{remark}

\begin{lemma}
$\G(L,\tilde{L})$ is well-defined up to Legendrian isotopy.
\end{lemma}

\begin{figure}
\centering{
\includegraphics[width=1in,angle=270]{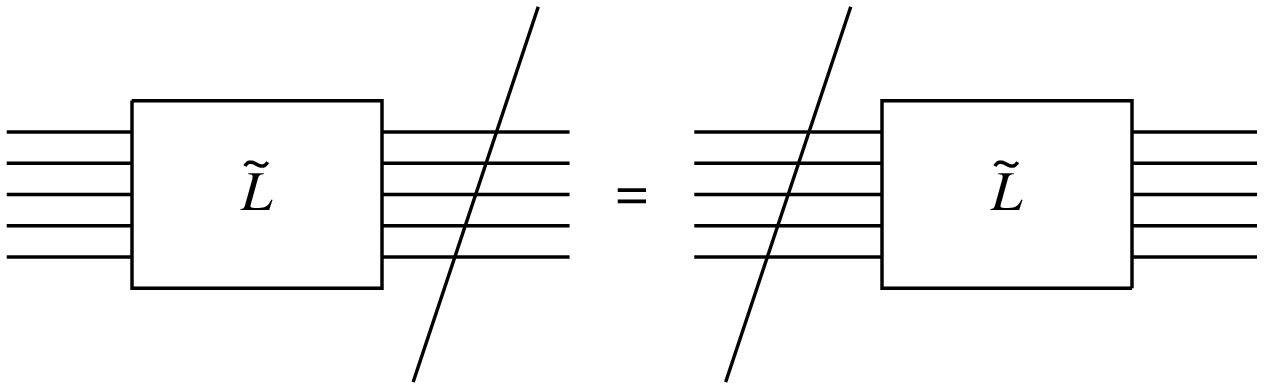}

\vspace{5pt}

\includegraphics[width=0.9in,angle=270]{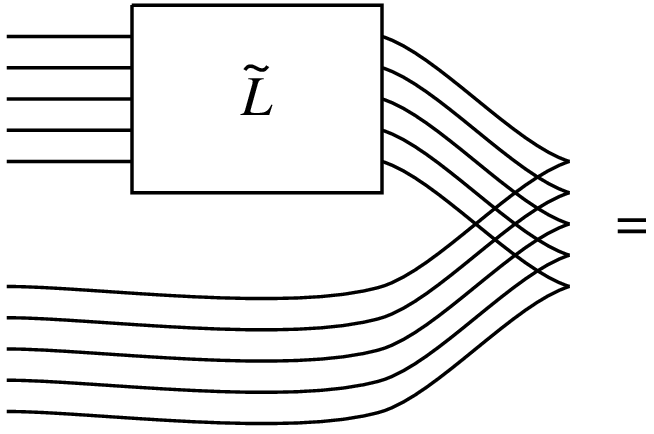}
\hspace{2pt}
\reflectbox{\includegraphics[width=0.9in,angle=270]{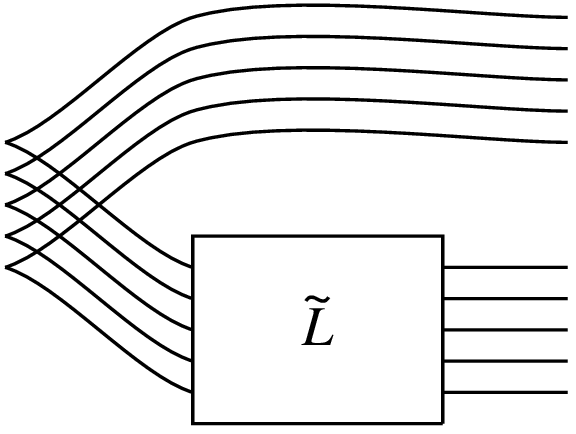}}
\hspace{12pt}
\includegraphics[width=0.9in,angle=270]{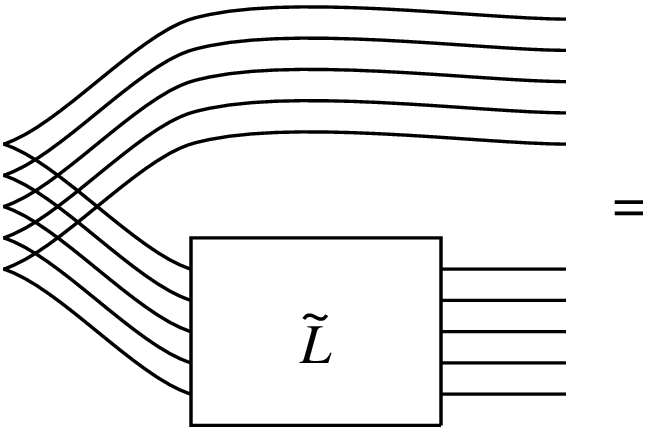}
\hspace{2pt}
\reflectbox{\includegraphics[width=0.9in,angle=270]{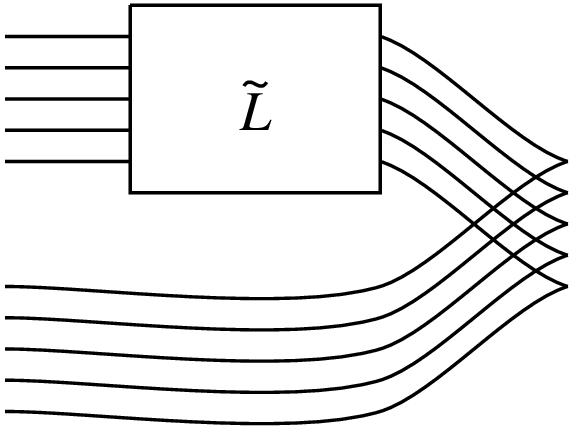}}
}
\caption{
Pushing $\tilde{L}$ through singularities in $L$: a crossing,
a right cusp, and a left cusp.
}
\label{fig:torusproof1}
\end{figure}

\begin{figure}
\centering{
\includegraphics[width=0.8in,angle=270]{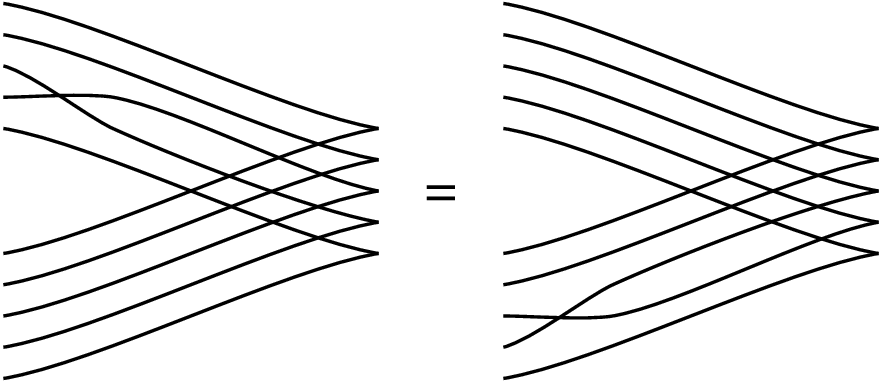}
\hspace{0.3in}
\includegraphics[width=0.8in,angle=270]{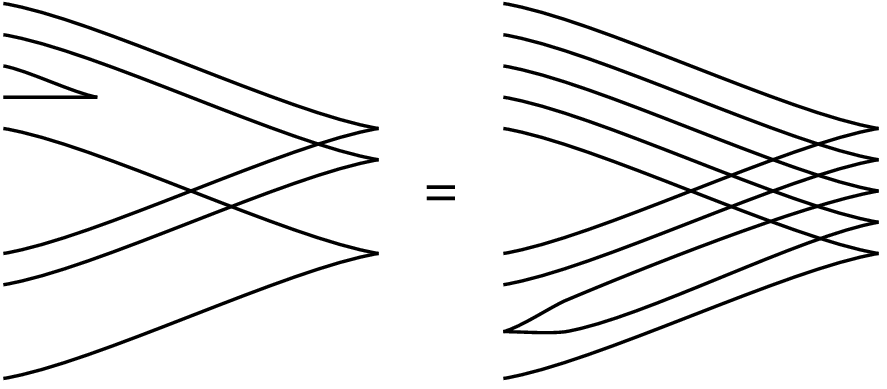}
\hspace{0.3in}
\includegraphics[width=0.8in,angle=270]{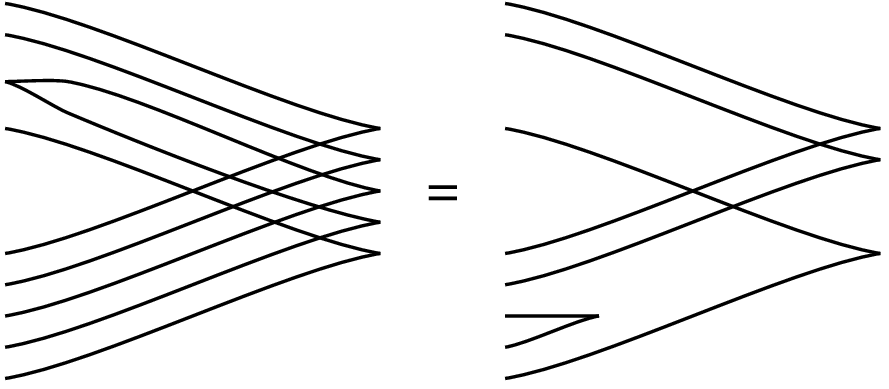}
}
\caption{
Pushing singularities in $\tilde{L}$ through a right cusp.
}
\label{fig:torusproof2}
\end{figure}

\begin{proof}
We need to show that, up to Legendrian isotopy, $\G(L,\tilde{L})$ is
independent of the piece of the $n$-copy of $L_1$ which we excise
and replace by $\tilde{L}$, as long as this piece is oriented
left to right.  The singularities of $\tilde{L}$ consist of
crossings, left cusps, and right cusps; we imagine pushing these
singularities one by one from one section of the $n$-copy of $L_1$ to
another.

We can clearly push these singularities through any piece of 
$\G(L,\tilde{L})$ which crosses a neighborhood of $\tilde{L}$ transversely;
see the top diagram in Figure~\ref{fig:torusproof1}.
Figure~\ref{fig:torusproof2} shows that we can also push 
singularities through a right cusp in $L_1$, and clearly this argument
extends to left cusps as well.  We conclude that we can push all
of $\tilde{L}$ through a cusp, resulting in the left-to-right
mirror reflection of $\tilde{L}$; see the bottom diagrams in
Figure~\ref{fig:torusproof1}.  The lemma follows.
\end{proof}

\begin{proposition}
$\G(L,\tilde{L})$ is a well-defined operation on Legendrian isotopy
classes; that is, if we change $L,\tilde{L}$ by Legendrian isotopies,
then $\G(L,\tilde{L})$ changes by a Legendrian isotopy as well.
\label{prop:torusgluing}
\end{proposition}

\begin{proof}
We first consider Legendrian-isotopy changes of $\tilde{L}$.
These fall into two categories: isotopies where the endpoints of
$\tilde{L}$ remain fixed, and horizontal translations of $\tilde{L}$
(i.e., moving the dashed lines).  The first category clearly preserves
the Legendrian isotopy class of $\G(L,\tilde{L})$.  The second category
consists of pushing singularities in $\tilde{L}$ through the dashed
lines.  But Figure~\ref{fig:torusproof2} shows that we can
push individual singularities from one side of $\tilde{L}$ to
the other, by moving the singularity all the way around the $n$-copy
of $L_1$.
Hence Legendrian-isotopy changes of $\tilde{L}$ do not change
$\G(L,\tilde{L})$.

\begin{figure}
\centering{
\includegraphics[width=0.9in,angle=270]{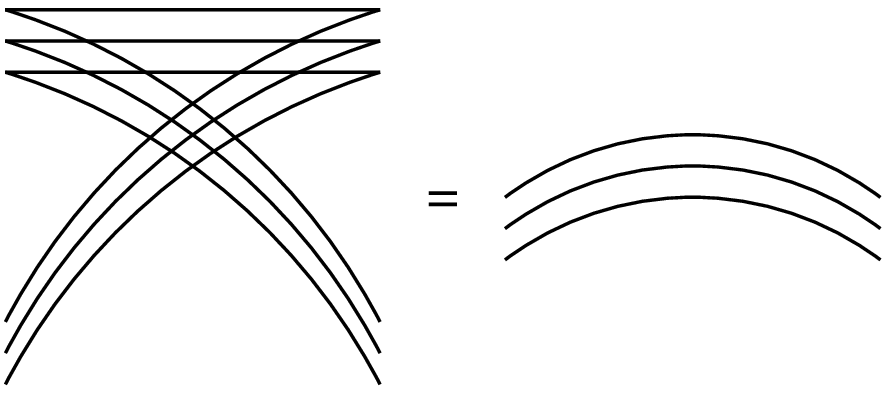}
\hspace{0.3in}
\includegraphics[width=0.9in,angle=270]{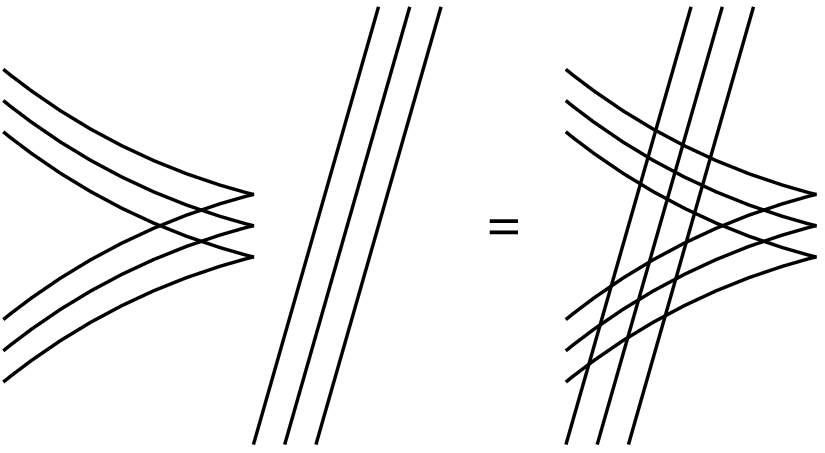}
\hspace{0.3in}
\includegraphics[width=0.9in,angle=270]{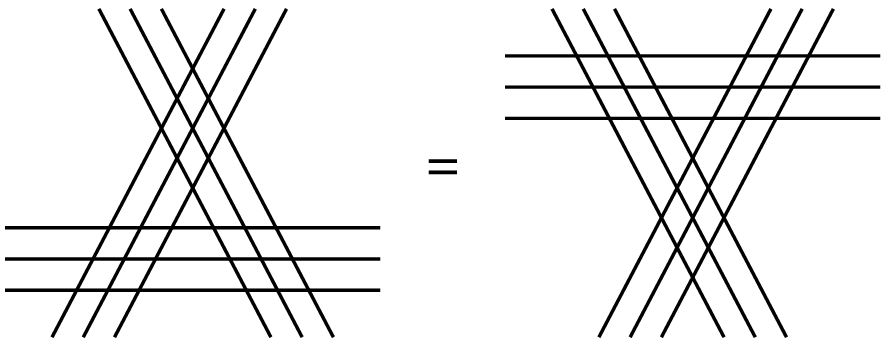}
}
\caption{
Legendrian Reidemeister moves on $n$-copies; in this illustration, $n=3$.
}
\label{fig:torusproof3}
\end{figure}

Next consider Legendrian-isotopy changes of $L$.  It suffices to
show that $\G(L,\tilde{L})$ does not change under Legendrian 
Reidemeister moves on $L$.  Consider such a move, and push $\tilde{L}$
away from a neighborhood of the move.  Then the fact that 
the Legendrian-isotopy class of $\G(L,\tilde{L})$ does not change
follows from Figure~\ref{fig:torusproof3}.
\end{proof}

\begin{remark}
Both Proposition~\ref{prop:torusgluing} and 
Corollary~\ref{cor:ncopy} below have been known for some time.
It is easy, and probably more natural, to establish
Proposition~\ref{prop:torusgluing} using the global, non-front
definition of Legendrian satellites; we chose to present
the front proof because of its concreteness.
\end{remark}

\begin{corollary}
Legendrian-isotopic knots in $\R^3$ have Legendrian-isotopic
$n$-copies.
\label{cor:ncopy}
\end{corollary}

\begin{figure}
\centering{
\includegraphics[height=5in,angle=270]{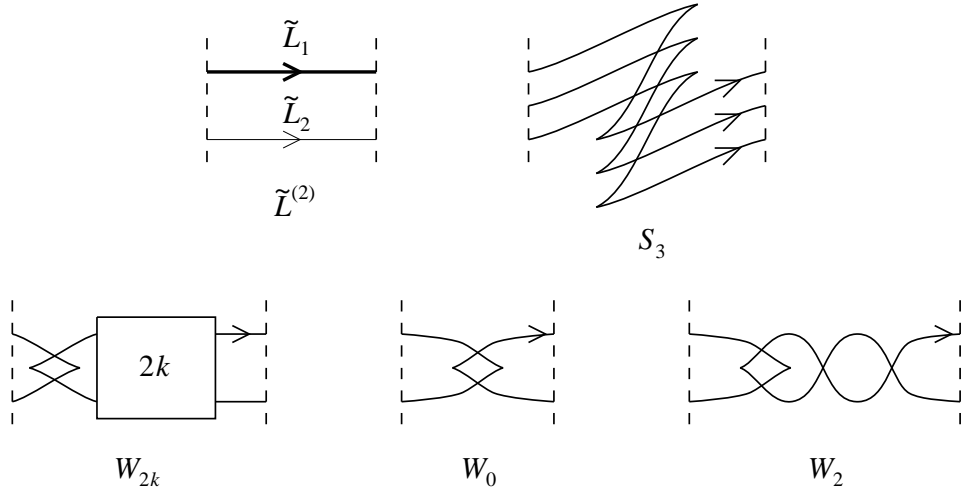}
}
\caption{
The solid-torus links $\tilde{L}^{(2)}$ and $S_3$ (with obvious
generalization to a family
of links $S_n$), and the solid-torus Whitehead knots
$W_{2k}$, $k \geq 0$, with $W_0$ and $W_2$ shown as examples.  The box
indicates $2k$ half-twists.
}
\label{fig:torusex}
\end{figure}

\begin{proof}
The $n$-copy of a knot $K$ is simply $\G(K,\tilde{L}^{(n)})$, where
$\tilde{L}^{(n)}$ is the union of $n$ unlinked loops which wind
once around $S^1 \times \R^2$; see Figure~\ref{fig:torusex}
for an illustration of $\tilde{L}^{(2)}$.
The result follows from Proposition~\ref{prop:torusgluing}.
\end{proof}

\begin{corollary}[\cite{bib:Mi}]
Suppose that $K$ is a stabilization of a Legendrian knot.
The $n$-copy of $K$ is Legendrian isotopic
to the $n$-copy with components cyclically permuted.
More precisely, if $L_1,\ldots,L_n$ are the components of the $n$-copy
of $K$, with $L_i$ slightly higher than $L_{i+1}$ in $z$ coordinate,
then $(L_1,L_2,\ldots,L_n)$ is Legendrian isotopic to
$(L_{1+k},L_{2+k},\ldots,L_{n+k})$ for any $k$, where indices
are taken modulo $n$.
\end{corollary}

\begin{proof}
Suppose, without loss of generality, that $K=S_+(K')$ for a Legendrian
knot $K'$.  Then the $n$-copy of $K$ is the Legendrian satellite 
$\G(K',S_n)$, where $S_n$ is the solid-torus ``$n$-copy 
stabilization link'' depicted in Figure~\ref{fig:torusex}.  
It is easy to see that $S_n$ is Legendrian isotopic
to itself with components cyclically permuted; now apply
Proposition~\ref{prop:torusgluing}.
\end{proof}

We now present some applications of Proposition~\ref{prop:torusgluing}
to knots and links on the solid torus.  Consider the link
$\tilde{L}^{(2)}$ shown in 
Figure~\ref{fig:torusex}.  The following result, established
in \cite{bib:Tra} using generating functions, is also
proven in \cite{bib:NT} using the DGA for solid-torus links.
The proof we give is yet another one.

\begin{proposition}
Write $\tilde{L}^{(2)} = (\tilde{L}_1,\tilde{L}_2)$.  Then 
$(\tilde{L}_1,\tilde{L}_2)$ is not Legendrian isotopic to
$(\tilde{L}_2,\tilde{L}_1)$.
\end{proposition}

\begin{proof}
In \cite[Proposition~4.11]{bib:Ng}, the author proves that the double 
of the figure
eight knot is not Legendrian isotopic to the double with components
swapped.  
The result now follows from Proposition~\ref{prop:torusgluing}.
\end{proof}

Now consider the Whitehead knots $W_{2k}$ shown in Figure~\ref{fig:torusex}.
Each $W_{2k}$ has $r=w=0$ and is thus topologically isotopic to its
inverse (the same knot with the opposite orientation).  
By contrast, we can now show the following result.

\begin{proposition}
$W_{2k}$ is not Legendrian isotopic to its inverse.
\label{prop:whiteheadknot}
\end{proposition}

\begin{proof}
As usual, write $-W_{2k}$ for the inverse of $W_{2k}$, and let $L$ be the
double of the usual ``flying-saucer'' unknot in $\R^3$.  For $k=1$,
it is easy to check that $\G(L,W_2)$ is precisely the oriented
Whitehead link from \cite[\S4.5]{bib:Ng}, and that $\G(L,-W_2)$ is
the same link with one component reversed.  
Proposition~\ref{prop:torusgluing} and \cite[Proposition~4.12]{bib:Ng}
then imply that $W_2$ and $-W_2$ are not Legendrian isotopic.

A calculation similar to the one in the proof of
\cite[Proposition~4.12]{bib:Ng}, omitted here, shows that $\G(L,W_{2k})$
and $\G(L,-W_{2k})$ are not Legendrian isotopic for arbitrary $k \geq 0$.
The result follows.
\end{proof}

The solid-torus DGA from \cite{bib:NT} fails to distinguish between
$W_{2k}$ and its inverse.  Proposition~\ref{prop:whiteheadknot}
is thus a result about solid-torus knots whose only presently known
proof uses the Legendrian satellite construction.

\section{Doubles}
\label{sec:doubles}

The Chekanov-Eliashberg 
DGA invariant
vanishes for links which are stabilizations.  The Legendrian satellite
construction, however, seems to yield nontrivial nonclassical invariants
of {\it all} Legendrian links; see Remark~\ref{rmk:stabinvts} below.  
On the other hand, the main result of this section shows that some of 
the simplest Legendrian satellites of stabilizations do not contain any 
new information.

\begin{definition}
The {\it Legendrian Whitehead double} of a Legendrian knot $K$ in $\R^3$
is $\G(K,W_0)$, where $W_0$ is the knot shown in Figure~\ref{fig:torusex}.
More generally, if $\tilde{L}$ has two endpoints, then
we call $\G(K,\tilde{L})$ a {\it satellite double} of $K$.
\end{definition}

As mentioned in Section~\ref{sec:gluing},
the Legendrian Whitehead double was originally defined by Eliashberg,
with further study by Fuchs \cite{bib:Fuc}, who uses the
notation $\Gamma_{\rm dbl}(0,0)$ for our $\G(K,W_0)$.

\begin{remark} {\it Legendrian satellites and maximal $tb$.}
By Remark~\ref{rmk:satellitetbr}, the Legendrian Whitehead double 
of any Legendrian knot has Thurston-Bennequin number 1.  
As noted by J.\ Sabloff and the author, it is easy to show that
the Legendrian Whitehead double maximizes $tb$ in its topological class.
This follows from the fact that $g(\G(K,W_0)) = 1$, along
with Bennequin's inequality $tb(K) \leq 2g(K)-1$ 
\cite{bib:Ben}, where $g(K)$ is the (three-ball) genus of $K$.  
A similar argument shows
that the usual double of any Legendrian knot maximizes $tb$.

It is not true, however, that all satellite doubles maximize $tb$,
even when $\tilde{L}$ maximizes $tb$.
In particular, if $\tilde{L}$ has a half-twist 
\includegraphics[height=12pt]{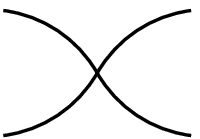}
next to its endpoints,
and $K$ is a stabilization, then $\G(K,\tilde{L})$ will also be a 
stabilization.
\end{remark}

\begin{proposition}
If $K_1$ and $K_2$ are stabilized Legendrian knots in the same
topological class with the same $tb$ and $r$, then
the DGAs of the Legendrian Whitehead doubles of $K_1$ and $K_2$
are equivalent.
\label{prop:whitehead}
\end{proposition}

The key to proving Proposition~\ref{prop:whitehead} is the following
result, whose proof we delay until Section~\ref{sec:prooflemdouble}.

\begin{lemma}
For any Legendrian knot $K$ which is a stabilization, the DGAs of 
$\G(K,W_0)$ and
of $\G(S_+S_-(K),W_0)$ are equivalent.
\label{lem:doubleproof}
\end{lemma}

\begin{proof}[Proof of Proposition \ref{prop:whitehead}]
By a result of \cite{bib:FT}, any two Legendrian knots which are
topologically identical and have the same $tb$ and $r$ are Legendrian
isotopic after some number of applications of the double stabilization
operator $S_+S_-$.  That is, there exists an $n \geq 0$ such that
$(S_+S_-)^n K_1$ and $(S_+S_-)^n K_2$ are Legendrian isotopic.
The proposition now follows directly from Lemma~\ref{lem:doubleproof}.
\end{proof}

\begin{remark}
It can in fact be shown that the DGA of the Legendrian Whitehead double
of a stabilized knot depends only on the $tb$ and $r$ of the knot,
and not on its topological class.  In particular, we can recover
the result of \cite{bib:Fuc} that the DGA of a Legendrian Whitehead double
always possesses an augmentation.
\end{remark}

A slightly modified version of the proof of Lemma~\ref{lem:doubleproof},
omitted here for simplicity,
establishes the following more general result.

\begin{proposition}
If $K_1$ and $K_2$ are stabilized Legendrian knots in the same
topological class with the same $tb$ and $r$, and $\tilde{L}$
is any Legendrian link in $S^1\times\R^2$ with two endpoints and
winding number zero, then the DGAs of $\G(K_1,\tilde{L})$
and $\G(K_2,\tilde{L})$ are equivalent.
\label{prop:torusdouble}
\end{proposition}

We believe that Proposition~\ref{prop:torusdouble} actually holds
for {\it any} satellite doubles of stabilized knots $K_1$ and $K_2$
with the same $tb$ and $r$, regardless of the 
winding number of $\tilde{L}$.  However, the analogue of
Lemma~\ref{lem:doubleproof} is false if $\tilde{L}$ has winding number
$\pm 2$, since $\G(K,\tilde{L})$ and $\G(S_+S_-(K),\tilde{L})$ have
different $tb$; see Remark~\ref{rmk:satellitetbr}.
Nevertheless, the argument of the proof of Lemma~\ref{lem:doubleproof}
shows that the characteristic algebra, at least, can never distinguish
between satellite doubles of stabilized knots.

\begin{remark} {\it Invariants of stabilized Legendrian knots.}
As mentioned in the Introduction, it remains a very interesting open problem to
find nonclassical invariants of stabilized Legendrian knots.  
There are currently
no methods to prove that two stabilized knots with the same topological
type, $tb$, and $r$ are not Legendrian isotopic.

We are hopeful that satellites more complicated than doubles will encode
interesting information for stabilized knots.
In particular, it seems that the $n$-copy of
any Legendrian link maximizes Thurston-Bennequin number when $n\geq 2$,
and thus probably has a nontrivial DGA.
By Corollary~\ref{cor:ncopy}, the DGAs of Legendrian satellites of a 
Legendrian link, including
the $n$-copy, are Legendrian-isotopy invariants, which likely contain
interesting nonclassical information in general.
The problem we face when dealing with complicated satellites, however, 
is extracting useful information from the Poincar\'e polynomials
or the characteristic algebra.  See \cite{bib:Mi}.

There is another approach to finding invariants
of stabilized knots, which is probably more natural than investigating
satellites.  Eliashberg, Givental, and Hofer \cite{bib:EGH} have recently
developed symplectic field theory, which generalizes contact homology;
\cite[\S2.8]{bib:EGH} describes how this method yields invariants
of Legendrian links, which would likely not vanish for stabilized knots.  
Unfortunately, no explicit combinatorial description, \`a la Chekanov, is 
presently known for the symplectic field theory associated to 
Legendrian links.
\label{rmk:stabinvts}
\end{remark}

\section{Proof of Lemma~\ref{lem:doubleproof}}
\label{sec:prooflemdouble}

\begin{figure}
\centering{
\includegraphics[height=4.5in]{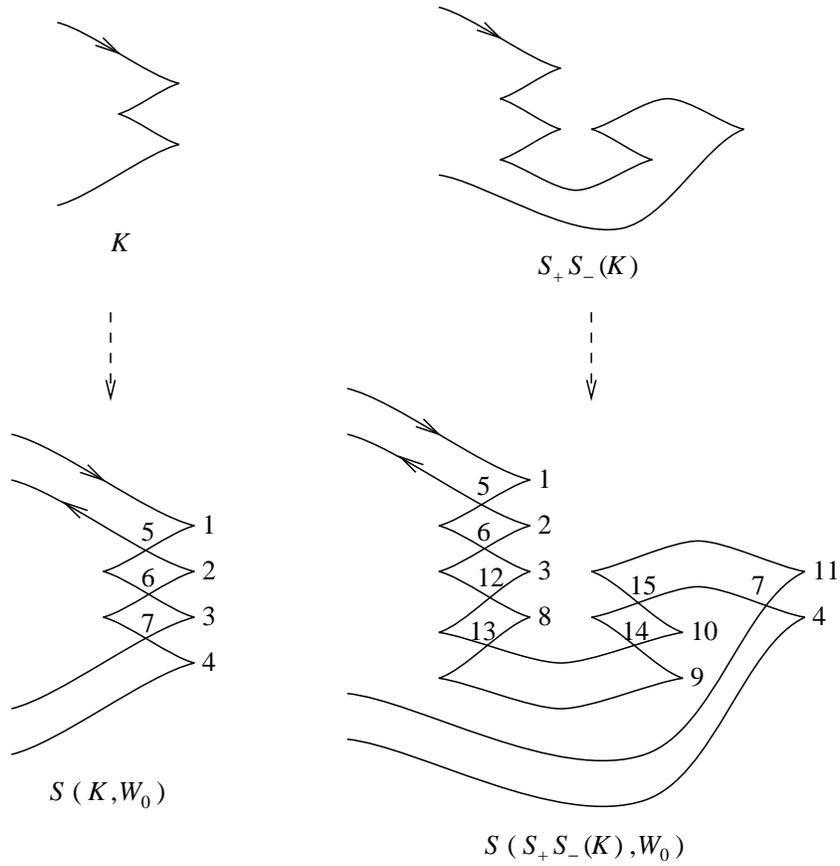}
}
\caption{
Whitehead doubles of stabilizations.  In the lower diagrams, 
vertices (crossings and right cusps) are labelled, with vertex $a_i$
labelled by $i$.
}
\label{fig:doubleproof}
\end{figure}

We assume familiarity with the definition of the DGA invariant,
as formulated in \cite{bib:Ng}.
We may suppose, without loss of generality, that $K = S_+(K')$ for some
Legendrian $K'$.  By using, if necessary, Legendrian Reidemeister 
moves (more precisely, IIb and the mirror of I from
\cite[Figure~1]{bib:Ng}), we may further
assume that the rightmost cusp in $K'$ is oriented downwards.
If we shift the zigzag in $K = S_+(K')$ next to the rightmost
cusp in $K'$, then $K$ and $\G(K,W_0)$ look like the diagrams
in Figure~\ref{fig:doubleproof} near the rightmost cusp.

The corresponding parts of $S_+S_-(K)$ and $\G(S_+S_-(K),W_0)$
are also shown in Figure~\ref{fig:doubleproof}, and 
$\G(K,W_0)$ and $\G(S_+S_-(K),W_0)$ are identical outside the
regions depicted.  It is easy to check that the degrees of all
vertices not depicted are equal for the two Legendrian Whitehead doubles,
and that the degrees of the vertices depicted are
$1$ for $a_1,a_2,a_3,a_4,a_8,a_9,a_{10},a_{11}$ and
$0$ for $a_5,a_6,a_7,a_{12},a_{13},a_{14},a_{15}$, in either
diagram.  Since the regions drawn are the rightmost parts of
each double, the DGA for $\G(S_+S_-(K),W_0)$
is simply obtained from the DGA for $\G(K,W_0)$ by making the following 
replacements:

\[
\left\{
\begin{array}{c}
\d a_2 = 1 - ta_5a_6 \\
\d a_3 = t^{-1} - a_6a_7 \\
\d a_5 = \d a_6 = \d a_7 = 0
\end{array}
\right\}
\qquad 
\longrightarrow
\qquad
\left\{
\begin{array}{c}
\d a_2 = 1 - ta_5a_6 \\
\d a_3 = t^{-1} - a_6a_{12} \\
\d a_8 = 1 - ta_{12}a_{13} \\
\d a_9 = t^{-1} - a_{14}a_{13} \\
\d a_{10} = 1 - ta_{15}a_{14} \\
\d a_{11} = t^{-1} - a_{15}a_7 \\
\d a_5 = \d a_6 = \d a_7 = 0 \\
\d a_{12} = \d a_{13} = \d a_{14} = \d a_{15} = 0
\end{array}
\right\} .
\]
We further note that none of the vertices depicted in
Figure~\ref{fig:doubleproof}, besides $a_1,a_4,a_5,a_7$, appears anywhere
in the DGAs except in the equations above; $a_5,a_7$ appear additionally
in $\d a_1, \d a_4$, respectively.

Our goal is to apply elementary automorphisms and algebraic stabilizations 
(see \cite{bib:Che}) to
the DGA for $\G(S_+S_-(K),W_0)$, until we obtain the DGA for
$\G(K,W_0)$.  Start with the DGA for $\G(S_+S_-(K),W_0)$;
we begin by rewriting $\d a_3, \d a_8, \d a_9, \d a_{10},
\d a_{11}$ in a more manageable form.

We first wish to rewrite $\d a_3$ as
$\d a_3 = t^{-1} - a_6 a_7$.  (Intuitively, this follows from the fact that
$a_5 = a_7$ in the characteristic algebra or in the homology of the DGA.)
We define the words $\alpha_1,\alpha_2,\alpha_3$ in the DGA as follows,
and then compute $\d \alpha_1, \d \alpha_2, \d \alpha_3$:
\[
\begin{aligned}
\alpha_1 &= ta_{15}a_9 - a_{10}a_{13} \\
\alpha_2 &= ta_{11} - \alpha_1a_7 \\
\alpha_3 &= a_8a_7 - a_{12}\alpha_2
\end{aligned}
\qquad \Longrightarrow \qquad
\begin{aligned}
\d \alpha_1 &= a_{15} - a_{13} \\
\d \alpha_2 &= 1 - t a_{13} a_7 \\
\d \alpha_3 &= a_7 - a_{12}.
\end{aligned}
\]
If we apply the elementary automorphism $a_3 \mapsto a_3 + a_6\alpha_3$,
then we obtain $\d a_3 = t^{-1} - a_6 a_7$.

In a similar fashion, we can successively replace
$\d a_{11}, \d a_{10}, \d a_9, \d a_8$ as follows:
$\d a_{11} = t^{-1} - a_{15} a_5$; $\d a_{10} = 1 - t a_6 a_{14}$;
$\d a_9 = t^{-1} - a_5 a_{13}$; $\d a_8 = 1 - t a_{12} a_6$.

For convenience, we now define
\[
\tilde{a}_2 = (1-ta_6a_5)a_3 + a_6a_2a_7
\qquad \Longrightarrow \qquad 
\d \tilde{a}_2 = t^{-1} - a_6a_5.
\]
By applying the elementary automorphisms
$a_{11} \mapsto a_{11} + \tilde{a}_2$ and $a_{15} \mapsto a_{15} + a_6$,
we obtain $\d a_{11} = - a_{15} a_5$.  Similarly, we may write
$\d a_{10} = -t a_6 a_{14}$, $\d a_9 = -a_5 a_{13}$, $\d a_8 = -t a_{12} a_6$.

At this point, the DGA has the following form:
\[
\begin{array}{rclcrcl}
\d a_2 &=& 1 - ta_5a_6 & \qquad \qquad &
\d a_9 &=& - a_5a_{13} \\
\d a_3 &=& t^{-1} - a_6a_7 &&
\d a_{10} &=& -ta_6a_{14} \\
\d a_8 &=& -ta_{12}a_6 &&
\d a_{11} &=& -a_{15}a_5 \\
\multicolumn{7}{c}{
\d a_5 = \d a_6 = \d a_7 = \d a_{12} = \d a_{13} = \d a_{14} = \d a_{15} = 0.
}
\end{array}
\]

We next eliminate $a_8,a_{11},a_{12},a_{15}$ through algebraic stabilization
and destabilization.
Introduce $e_1$ and $e_2$ of degree $0$ and $-1$, respectively, with
$\d e_1 = e_2$ (and $\d e_2 = 0$).  Let $\Phi_1$ be the composition
of the following elementary automorphisms in succession:
\[
a_{12} \mapsto a_{12} - e_1a_5; \qquad
a_8 \mapsto a_8 - e_1 a_2; \qquad
e_1 \mapsto e_1 + t a_{12}a_6 - e_2a_2.
\]
Under $\Phi_1$, the DGA changes as follows:
\[
\left\{
\begin{aligned}
\d a_8 &= -ta_{12} a_6 \\
\d e_1 &= e_2 \\
\d a_{12} &= 0 \\
\d e_2 &= 0
\end{aligned}
\right\}
\qquad \stackrel{\Phi_1}{\longrightarrow} \qquad
\left\{
\begin{aligned}
\d a_8 &= e_1 \\
\d e_1 &= 0 \\
\d a_{12} &= e_2a_5 \\
\d e_2 &= 0
\end{aligned}
\right\} .
\]
We may then drop $a_8$ and $e_1$; these simply correspond to an
algebraic stabilization.

Let $\Phi_2$ be the composition of the following maps:
\[
a_{15} \mapsto a_{15} - t a_{12}a_6; \quad \negthickspace\negthickspace
a_{11} \mapsto a_{11} - t a_{12} \tilde{a}_2; 
\quad \negthickspace\negthickspace
a_{12} \mapsto a_{12} + a_{15}a_5 - te_2a_5\tilde{a}_2; 
\quad \negthickspace\negthickspace
a_{15} \mapsto a_{15} + e_2a_2.
\]
Under $\Phi_2$, the DGA now changes as follows:
\[
\left\{
\begin{aligned}
\d a_{11} &= a_{15}a_5 \\
\d a_{12} &= e_2a_5 \\
\d a_{15} &= 0 \\
\d e_2 &= 0
\end{aligned}
\right\}
\qquad 
\stackrel{\Phi_2}{\longrightarrow} \qquad
\left\{
\begin{aligned}
\d a_{11} &= a_{12} \\
\d a_{12} &= 0 \\
\d a_{15} &= e_2 \\
\d e_2 &= 0
\end{aligned}
\right\} .
\]
We can now drop $a_{11},a_{12},a_{15},e_2$; these correspond to two
algebraic stabilizations.

Hence, up to algebraic stabilizations, we have eliminated
$a_8,a_{11},a_{12},a_{15}$.  An entirely similar process allows us to
eliminate $a_9,a_{10},a_{13},a_{14}$.  The resulting DGA is
precisely the DGA of $\G(K,W_0)$, as desired. \hfill $\Box$

\vspace{12pt}

\begin{remark}
The only part of this proof which uses the structure of
$W_0$ is the calculation of the degrees of the vertices in
Figure~\ref{fig:doubleproof}.  To prove the more general case given in
Proposition~\ref{prop:torusdouble}, we have to take more care
vis-\`a-vis degrees, but the idea is the same.
The proof also extends to knots which are not satellite doubles,
but whose rightmost parts look like the bottom diagrams in
Figure~\ref{fig:doubleproof}.
\end{remark}

\begin{remark}
The method of the proof can also be used to show that the invariant
$HC_{132}$ introduced by Michatchev \cite{bib:Mi} does not encode any
nonclassical information for stabilized knots.  The computation for 
this case involves more generators, but less algebra, than the 
computation performed in this section.
\end{remark}

\vspace{22pt}


\vspace{24pt}

\end{document}